\documentclass[graybox]{svmult}

\usepackage{mathptmx}       
\usepackage{helvet}         
\usepackage{courier}        
\usepackage{type1cm}        
\usepackage{makeidx}         
\usepackage{graphicx}        
\usepackage{multicol}        
\usepackage[bottom]{footmisc}

\usepackage{amsmath, amssymb, latexsym}
\usepackage{graphicx}
\usepackage{enumerate}
\usepackage{psfrag}
\usepackage{array}
\usepackage{graphicx}
  \usepackage{pgfplots}
  \pgfplotsset{compat=newest}
  \usetikzlibrary{plotmarks}
\usetikzlibrary{matrix,shapes,arrows,positioning,fit,chains,scopes}
  \usepackage{grffile}
  \usepackage{mathtools}

\newcommand{\figureheight}{5.5cm}
\newcommand{\figurewidth}{8cm}
\newcommand{\yticklabelstyle}{yticklabel style={text width=0.65em,align=right}}
\newcommand{\xticklabelstyle}{} 

\definecolor{darkpastelgreen}{rgb}{0.01, 0.75, 0.24}
\definecolor{darkgreen}{rgb}{0.0, 0.2, 0.13}
\definecolor{ao}{rgb}{0.0, 0.5, 0.0}
\definecolor{ay}{rgb}{0.7, 0.7, 0.0}
\definecolor{al}{rgb}{0.8, 0.4, 0.0}

\definecolor{blue}{rgb}{0.8, 0.4, 0.0}

\newcommand{\slfrac}[2]{{\left.#1\middle/#2\right.}}

\newcommand{\R}{\mathbb R}

\newcommand{\dt}{{\Delta t}}
\newcommand{\dx}{{\Delta x}}

\DeclareMathOperator{\RS}{RS}

\DeclareMathOperator{\KT}{KT}
\DeclareMathOperator{\BE}{BE}
\DeclareMathOperator{\CN}{CN}

\DeclareMathOperator{\mm}{minmod}
\DeclareMathOperator{\sgn}{sgn}

\makeindex             

\begin{document}

\title*{A splitting approach for freezing waves.}
\author{Robin Flohr and Jens Rottmann-Matthes}
\institute{Robin Flohr \at Institute for Analysis, Karlsruhe Institute of Technology\\ Englerstra\ss e 2, 76131 Karlsruhe, Germany\\ \email{Robin.Flohr@kit.edu}
\and Jens Rottmann-Matthes \at Institute for Analysis, Karlsruhe Institute of Technology\\ Englerstra\ss e 2, 76131 Karlsruhe, Germany\\ \email{ Jens.Rottmann-Matthes@kit.edu }}
%
%
\maketitle

\abstract*{We present a numerical method which is able to approximate
  traveling waves (e.g.\ viscous profiles) in systems with hyperbolic
  and parabolic parts by a direct long-time forward
  simulation.  A difficulty with long-time simulations of traveling
  waves is that the solution leaves any bounded computational domain in
  finite time. To
  handle this problem one should go into a suitable co-moving frame.
  Since the velocity of the wave is typically unknown, we use the method of
  freezing \cite{BT04}, see also \cite{beynOttenRM2014},
  which transforms the original partial
  differential equation (PDE) into a partial differential algebraic
  equation (PDAE) and calculates a suitable co-moving frame
  on the fly.  The efficient numerical approximation of this freezing
  PDAE is a challenging problem and we introduce a new
  numerical discretization, suitable for problems that consist of
  hyperbolic conservation laws which are (nonlinearly) coupled to
  parabolic equations.
  The idea of our scheme is to use the operator splitting approach.
  The benefit of splitting methods in this context lies in the
  possibility to solve hyperbolic and parabolic parts with different
  numerical algorithms.\\
  We test our method at the (viscous) Burgers'
  equation.  Numerical experiments show linear and
  quadratic convergence rates for the approximation of the numerical
  steady state obtained by a long-time simulation 
  for Lie- and Strang-splitting respectively.  Due to these affirmative
  results we expect our scheme to be suitable also for freezing waves in
  hyperbolic-parabolic coupled PDEs.}

\abstract{We present a numerical method which is able to approximate
  traveling waves (e.g.\ viscous profiles) in systems with hyperbolic
  and parabolic parts by a direct long-time forward
  simulation.  A difficulty with long-time simulations of traveling
  waves is that the solution leaves any bounded computational domain in
  finite time. To
  handle this problem one should go into a suitable co-moving frame.
  Since the velocity of the wave is typically unknown, we use the method of
  freezing \cite{BT04}, see also \cite{beynOttenRM2014},
  which transforms the original partial
  differential equation (PDE) into a partial differential algebraic
  equation (PDAE) and calculates a suitable co-moving frame
  on the fly.  The efficient numerical approximation of this freezing
  PDAE is a challenging problem and we introduce a new
  numerical discretization, suitable for problems that consist of
  hyperbolic conservation laws which are (nonlinearly) coupled to
  parabolic equations.
  The idea of our scheme is to use the operator splitting approach.
  The benefit of splitting methods in this context lies in the
  possibility to solve hyperbolic and parabolic parts with different
  numerical algorithms.\\
  We test our method at the (viscous) Burgers'
  equation.  Numerical experiments show linear and
  quadratic convergence rates for the approximation of the numerical
  steady state obtained by a long-time simulation 
  for Lie- and Strang-splitting respectively.  Due to these affirmative
  results we expect our scheme to be suitable also for freezing waves in
  hyperbolic-parabolic coupled PDEs.}

\section{Introduction}
Many partial differential equations from applications consist of
different parts,
\begin{equation}\label{flohr:0}u_t=Au_{xx}+f(u)_x+g(u) \eqqcolon F(u).
\end{equation}
Sometimes, one part is parabolic while another part is hyperbolic
and these parts are nonlinearly coupled.  Examples of such
hyperbolic-parabolic coupled PDEs are hyperbolic models of
chemosensitive movement or reaction-diffusion equations for which not
all components diffuse.  

One is often interested in special solutions, which
arise as (time-)asymptotic limits of solutions to the
Cauchy problem for \eqref{flohr:0}.  An important class of such solutions
are traveling waves.  They describe how mass (or information) travels
through the domain.  From this interpretation, it is obvious,
that one is often interested not only in the shape but also the velocity
of the traveling wave.

Traveling waves are solutions of the PDE \eqref{flohr:0} of the form
\begin{equation*}
u(x,t) = \bar v(x- \bar \mu t), \quad x\in\R, \quad t\in\R,
\end{equation*}
where $\bar v:\R \to \R$ is the non-constant profile and $\bar \mu\in\R$
the velocity of the wave.
For Burgers' equation there is a family of traveling wave solutions,
\begin{equation}
\label{flohr:twSolutionBurgers}
\begin{aligned}
u(x,t) &= \varphi(x-\bar\mu t) + \tfrac{1}{2} (b+c)=\bar v(x-\bar\mu t),\\
 \qquad \varphi(x) &= a\frac{1 - e^{ax}}{1 + e^{ax}},\qquad a=\tfrac{1}{2}(b-c), \qquad \bar\mu=\tfrac{1}{2}(b+c),
\end{aligned}
\end{equation}
parametrized by the asymptotic states
$\lim_{x\to-\infty} u(x,t) = b > c = \lim_{x\to\infty} u(x,t)$. 

As a toy example we consider the Cauchy problem for the
viscous Burgers' equation
\begin{equation}
\label{flohr:burgersequation}
u_t + (\tfrac{1}{2}u^2)_x = u_{xx}\; \text{ on }
\R\times[0,\infty),\quad
u(\cdot,0) = u_0\; \text{ on } \R.
\end{equation}

For the numerical approximation of \eqref{flohr:burgersequation} 
one has to truncate the unbounded spatial domain to a finite interval.
This leads to the problem that every traveling wave solution with
non-zero speed eventually leaves the computational domain.
The simplest remedy is to use periodic boundary conditions on a very large domain, but this is
only reasonable in the case of pulses.
Instead, we use the freezing method from \cite{BT04}.  The idea is to
move the spatial frame with the speed of the traveling wave.
We make the \textbf{ansatz} that the solution of \eqref{flohr:0} is of
the form
\begin{equation}\label{flohr:ansatz}
  \begin{aligned}
    u(x,t) &= v\big(x-{\color{black}\gamma(t)},t\big), \qquad {\color{black}\gamma(t)} \in \R,
  \end{aligned}
\end{equation}
where $\gamma(t)$ denotes a time dependent position. Then,
$\mu(t):=\gamma_t(t)$ can be interpreted as the velocity of the wave at
time $t$.
Plugging \eqref{flohr:ansatz} into \eqref{flohr:0} yields
\begin{equation}
\label{flohr:freezingeq}
v_t = F(v) + \mu(t)v_x,
\end{equation}
where both $v$ and $\mu$ are unknown.
Due to the additional unknown $\mu$ one has to
complement \eqref{flohr:freezingeq} with an addition algebraic equation,
called phase condition in \cite{BT04}, to retain well-posedness.
In Burgers' case this transforms \eqref{flohr:burgersequation} into the PDAE
\begin{equation}
\label{flohr:freezedBurgers}
\left\{
\begin{aligned}
v_t &= v_{xx} - (\tfrac{1}{2}v^2)_x + \mu v_x,\\
0 &= \Psi(v,\mu),\\
\gamma_t &= \mu(t),
\end{aligned}\right. \qquad
\begin{aligned}
v(\cdot,0) &= u_0,\\
 \gamma(0)&=0.
\end{aligned}
\end{equation}
This is called the \emph{freezing method} in \cite{BT04}.
We restrict to two standard choices for the phase condition,
the \emph{orthogonal phase condition} given by
\begin{equation}
\label{pco}
0 = \Psi(v,\mu) \coloneqq \langle v_t \;|\; v_x\rangle_{L^2} = \langle v_{xx} -(\tfrac12 v^2)_x + \mu v_x \;|\; v_x\rangle_{L^2}
\end{equation}
and the \emph{fixed phase condition} given by
\begin{equation}
\label{pcf}
0 = \Psi(v,\mu) \coloneqq \langle v-\hat v \;|\; \hat v_x\rangle_{L^2}
\end{equation}
with $\hat v$ an appropriately chosen reference function.

For the numerical approximation of \eqref{flohr:freezedBurgers} we use
splitting methods, which we briefly recall for convenience.
Assume that a solution to an initial value
problem of the form
\begin{equation}
\label{flohr:eqtoSplit}
u_t = A(u) + B(u)
\end{equation}
 is sought.  Let 
$\Phi_A^t(u_0)$ and $\Phi_B^t(u_0)$ denote the solution operators for $u_t =
A(u)$ and $u_t = B(u)$ with initial value $u_0$, respectively.
The \emph{Lie-Trotter splitting},
\begin{equation}\label{flohr:LieSplitting}
u^{n+1} = {\color{black}\Phi_B^\dt} \circ {\color{black}\Phi_A^\dt}(u^n),
\end{equation}
typically converges linearly to the exact solution for $\dt
\to 0$.
A splitting method that typically leads to second order convergence is
\emph{Strang splitting},
\begin{equation}\label{flohr:StrangSplitting}
u^{n+1} = {\color{black}\Phi_A^\slfrac{\dt}{2}} \circ {\color{black}\Phi_B^\dt} \circ {\color{black}\Phi_A^\slfrac{\dt}{2}}(u^n).
\end{equation}
For example in \cite{HoldenLubichRisebro}, the authors show that this
scheme is second order convergent for the viscous Burgers' PDE.

To apply this approach to the freezing PDAE
\eqref{flohr:freezedBurgers}, we split the
equation into two parts to separate the hyperbolic and parabolic
problem.  Then we solve each part with a method which is
particularly adapted to the respective sub-problem.  Namely we solve the
hyperbolic problem with an explicit scheme from
Kurganov and Tadmor \cite{kurganovTadmor}.
The parabolic sub-problem is solved by an implicit second order
finite-difference approximation, due to the restrictive CFL condition.

Our main focus in this article is on approximating the limits of the
time evolution and, unlike \cite{HoldenLubichRisebro},
not on the finite time convergence properties of the scheme.  In
particular, we aim to understand the preservation of
steady states and their stability by our schemes. 

%
In the case of ordinary differential equations there is a
well-established theory for numerical steady states.
For example in \cite{stuart1996} there are results which state that 
one-step methods preserve fixed points and their stability in a
$\dt^r$-shrinking neighborhood under Lipschitz assumptions,
where $r$ denotes the order of consistency of the one-step method.
An analogous result holds for the Strang splitting:
\begin{theorem}[\cite{F13}]
Let $A,B \in C^3(\R^m,\R^m)$ and assume that $\hat u$ is a hyperbolic
fixed point of \eqref{flohr:eqtoSplit}.
Let $\varphi_A,\varphi_B$ be one-step methods approximating
$\Phi_A,\Phi_B$ respectively.  If $\varphi_A,\varphi_b$ are second order
Runge-Kutta methods then there exist $\dt_0,K > 0$, such that the
numerical Strang splitting, $U^{n+1} = \varphi^{\dt}(U^n)=\varphi^{\slfrac{\dt}{2}}_B \circ
\varphi^\dt_A \circ \varphi^{\slfrac{\dt}{2}}_B(U^n)$, has a fixed point
$\hat U$ which is unique in the ball $B(\hat u;K\dt^2)$ for all $0<\dt
\le \dt_0$.
Furthermore, $\hat U$ is a stable (resp.\ unstable) fixed point of
$\varphi^\dt$ if $\hat u$ is a stable (resp.\ unstable) steady state of
\eqref{flohr:eqtoSplit}.
\end{theorem}

For the freezing method there are several results available,
where it is shown that the (continuous) method provides good approximations
including preservation of asymptotic stability 
of traveling waves for certain problem classes in the continuous and
semi-discrete case, see \cite{jens2012,beynOttenRM2014,thuemmlerDiss}.
But the time-asymptotic behavior of a discretization with a splitting
approach was never considered before.

A different approach to apply adapted schemes for different
parts of the freezing PDAE appears in \cite{jens2016}, where the
freezing method is used to capture similarity solutions of the
multidimensional Burgers' equation.  There
an IMEX-Runge-Kutta approach is used and
second oder convergence for the time dependent problem is shown on
finite time intervals. 

\section{The splitting scheme}

We now explicitly state our numerical scheme.
We split \eqref{flohr:freezedBurgers} into two sub-problems as follows:
Let $\Phi_A^t : (z_0,\gamma_0,\mu_0) \mapsto \big( z(t),\gamma(t),\mu(t) \big)$ be the solution operator
to the parabolic problem
\begin{equation}
\label{flohr:A}
\tag{A}
\left\{
\begin{aligned}
z_t &= z_{xx},\\
\gamma_t &= 0,\\
\mu_t &= 0,
\end{aligned}\right. \qquad
\begin{aligned}
z(\cdot,0) &= z_0,\\
 \gamma(0)&=\gamma_0,\\
  \mu(0)&=\mu_0,
\end{aligned}
\end{equation}
let $\Phi_B^t : (w_0,\gamma_0,\mu_0) \mapsto \big( w(t),\gamma(t),\mu(t) \big)$
be the solution operator to
\begin{equation}
\label{flohr:B}
\tag{B}
\left\{
\begin{aligned}
w_t &= - (\tfrac{1}{2}w^2)_x + \mu w_x,\\
0 &= \Psi(w,\mu),\\
\gamma_t &= \mu(t),
\end{aligned}\right. \qquad
\begin{aligned}
w(\cdot,0) &= w_0,\\
 \gamma(0)&=\gamma_0.
\end{aligned}
\end{equation}
Here $\Psi$ is one of the phase conditions \eqref{pco} or \eqref{pcf}.
Note that the initial value $\mu_0$ is ignored for this operator \eqref{flohr:B},
because it is uniquely determined by the constraints.
Since the splitting approach now iterates both solution operators consecutively, 
the question when and how to solve the algebraic constraint arises.
For the orthogonal phase condition we chose an explicit and
for the fixed phase condition we use a half-explicit approach.
Thus we calculate the speed $\mu$ prior to solving the nonlinear PDE,
the $\mu w_x$ part is then discretized by using finite differences.
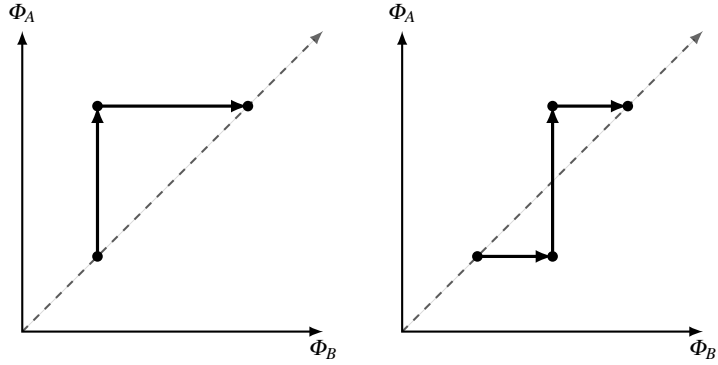
\begin{figure}[tb]
\centering

\begin{minipage}{1\linewidth}
\begin{minipage}{0.48\linewidth}
\hfill
\begin{tikzpicture}[thick,->,>=latex]
%
%
 \draw [<->,thick] (0,4) node (yaxis) [above] {$\Phi_A$}
         |- (4,0) node (xaxis) [below] {$\Phi_B$};
  
     \coordinate (a) at (1,1);
     \coordinate (b) at (1,3);
     \coordinate (c) at (3,3);
     
     \coordinate (A) at (0,0);
     \coordinate (B) at (4,4);
 
      \fill[black] (a) circle (2pt);
      \fill[black] (b) circle (2pt);
      \fill[black] (c) circle (2pt);
 
     \draw[->,solid,very thick] (a) -- (b);
     \draw[->,solid,very thick] (b) -- (c);

     \draw[dashed,opacity=0.6] (A) -- (B);

 \end{tikzpicture}
\end{minipage}
\hfill
\begin{minipage}{0.48\linewidth}
\begin{tikzpicture}[thick,->,>=latex]
 %
 %
  \draw [<->,thick] (0,4) node (yaxis) [above] {$\Phi_A$}
          |- (4,0) node (xaxis) [below] {$\Phi_B$};
   
      \coordinate (a) at (1,1);
      \coordinate (b) at (2,1);
      \coordinate (c) at (2,3);
      \coordinate (d) at (3,3);
      
      \coordinate (A) at (0,0);
      \coordinate (B) at (4,4);
  
       \fill[black] (a) circle (2pt);
       \fill[black] (b) circle (2pt);
       \fill[black] (c) circle (2pt);
       \fill (d) circle (2pt);
  
      \draw[->,solid,very thick] (a) -- (b);
      \draw[->,solid,very thick] (b) -- (c);
      \draw[->,solid,very thick] (c) -- (d);

      \draw[dashed,opacity=0.6] (A) -- (B);

  \end{tikzpicture}
\end{minipage}
\end{minipage}
 \caption{Diagram of the Lie-Trotter splitting on the left and the Strang splitting on the right.}
\label{fig:figuresSplittigDiagram}
\end{figure}
Lie and Strang splitting are illustrated by classical diagrams in \textbf{Fig.~\ref{fig:figuresSplittigDiagram}}.
A step in the vertical direction in Fig.~\ref{fig:figuresSplittigDiagram} amounts in numerically solving the Cauchy problem for the heat equation
\eqref{flohr:A}, whereas a step in the horizontal direction amounts to solve the hyperbolic PDAE \eqref{flohr:B}.
Only states on the dashed diagonal line might be considered as approximations to solutions to the original problem.
In addition, 
the order of the sub-problems \eqref{flohr:A}, \eqref{flohr:B} in the splitting approach is chosen such that the phase condition is satisfied at the end of a full time step.
More details about how to calculate the speed with the algebraic constraint can be found in the description of the schemes, see \eqref{flohr:phiBRSO}, \eqref{flohr:phiBRSF}, \eqref{flohr:phiBKTF} and \eqref{flohr:phiBKTO}.
A schematic overview of the schemes is given in Fig.~\ref{fig:table_used_numerical_schemes}.
\begin{figure}[tb]
\centering
{\setlength{\extrarowheight}{7pt}%
\begin{tabular}{| l | c | c | c | c |}
\hline
  {\setlength{\extrarowheight}{0pt} \begin{tabular}{@{}l@{}}\textbf{convergence}\end{tabular}} & \multicolumn{2}{c |}{\textbf{order 1}} & \multicolumn{2}{c|}{\textbf{order 2}}\\[7pt]
    \hline
  {\setlength{\extrarowheight}{0pt} \begin{tabular}{@{}l@{}}full problem\end{tabular}} & \multicolumn{4}{c |}{$u_t + \big(\tfrac{1}{2} u^2\big)_x = u_{xx}$}\\[7pt]
  \hline
  {\setlength{\extrarowheight}{0pt} \begin{tabular}{@{}l@{}}freezing method\end{tabular}} & \multicolumn{2}{c|}{orthogonal or fixed p.c.} & \multicolumn{2}{c|}{fixed p.c.}\\[7pt]
  \hline
 {\setlength{\extrarowheight}{0pt} \begin{tabular}{@{}l@{}}sub-problem\end{tabular}} & $w_t = -(\tfrac{1}{2}w^2)_x$ & $z_t =z_{xx}$ &  $w_t = -(\tfrac{1}{2}w^2)_x$ & $z_t =z_{xx}$\\[7pt]
   \hline
 {\setlength{\extrarowheight}{0pt} \begin{tabular}{@{}l@{}}semi-discrete \\ formulation\end{tabular}} & Rusanov Scheme & discrete Laplacian & Kurganov-Tadmor & discrete Laplacian\\[7pt]
    \hline
    {\setlength{\extrarowheight}{0pt} \begin{tabular}{@{}l@{}}time \\ discretization\end{tabular}} & forward Euler& backward Euler& Heun's method& Crank-Nicolson\\[7pt]
    \hline
       {\setlength{\extrarowheight}{0pt} \begin{tabular}{@{}l@{}}splitting method\end{tabular}} & \multicolumn{2}{c|}{Lie} & \multicolumn{2}{c|}{Strang}\\[7pt]
       \hline
\end{tabular}}
\caption{Overview of the applied numerical schemes for the presented schemes which offer a numerical steady state.}
\label{fig:table_used_numerical_schemes}
\end{figure}

\subsection{First order scheme.}

We first present a first order scheme. For this we use a method of lines (MOL) approach for \eqref{flohr:B}:
We choose a finite interval $[L_-,L_+]$ and choose a spatial grid with uniform step size $\dx$ and spatially discretize with the semi-discrete version of the Rusanov scheme.
We use Dirichlet boundary conditions.
%
It is worth mentioning here that for example the LxF or NT scheme do not offer a semi-discrete version.
The \textbf{Rusanov scheme (RS)} in its semi-discrete form is given by
\begin{equation*}
\begin{aligned}
\frac{d}{dt} w_j(t) &= - \frac{f\big(w_{j+1}(t)\big) - f\big(w_{j-1}(t)\big)}{2 \Delta x}\\
&\qquad + \frac{\kappa}{2 \Delta x} [w_{j+1}(t) - 2 w_j(t) + w_{j-1}(t)]\\
&= - \partial_0 f\big(w(t)\big)_j + \kappa \tfrac{\Delta x}{2} \partial_0^2 w(t)_j\\
&\eqqcolon \RS^\dx(w(t)),
\end{aligned}
\end{equation*}
where $\partial_0$ is the central difference quotient, $\partial_0 w_j = \tfrac{1}{2\dx}(w_{j+1} - w_{j-1})$,
$\partial_0^2$ the discrete Laplacian, both with Dirichlet boundary conditions and $\kappa = \max_j u(j\dx,0)$ is the maximum over the initial value evaluated at all grid points.

This scheme is in a simplified form: Since the local maximal speeds, used in the Kurganov-Tadmor scheme, ensure that all information of the Riemann fans stay in each cell of the discretized problem, they can be replaced by an upper bound. In the case of the Burgers' nonlinearity this upper bound is given by the maximal absolute value of the solution, which, in turn, is given by the maximal absolute value $\kappa$ of the initial function $u_0$ due to the maximums principle.

The time discretization is done with a uniform step size $\dt$, for the first order version we use the forward Euler method.
The numerical approximation of $\Phi_B^\dt$ will be denoted by
$\phi_{B,\RS\dx}^\dt$ and $\varphi_{B,\RS\dx}^\dt$ for the two different phase conditions \eqref{pco} and \eqref{pcf}, respectively.
The operator $\phi_{B,\RS\dx}^\dt$ is given as the function which maps $w_0,\gamma_0,\mu_0$ to the solution $w^1,\gamma^1,\mu^1$ of
the system
\begin{equation}
\label{flohr:phiBRSO}
\left\{
\begin{aligned}
w^{1} &= w^0 + \dt \RS^\dx (w^0) + \dt\mu^*\partial_0w^0,\\
\mu^* &= - \frac{\partial_0 w^{0\top} \big( \partial_0^2 w^0 - w^0\partial_0 w^0 \big)}{\partial_0w^{0\top}\partial_0w^0},\\
\gamma^1 &= \gamma_0 + \dt\mu^1,\\
\mu^1 &= - \frac{\partial_0 w^{1\top} \big( \partial_0^2 w^1 - w^1\partial_0 w^1 \big)}{\partial_0w^{1\top}\partial_0w^1},
\end{aligned}\right. \qquad
\begin{aligned}
w^0 &= w_0,\\
 \gamma^0&=\gamma_0,
\end{aligned}
\end{equation}
where we use a discrete version of the orthogonal phase condition \eqref{pco}.
For the fixed phase condition  \eqref{pcf} the operator $\varphi_{B,\RS\dx}^\dt$ is given
as the mapping, which maps  $w_0,\gamma_0,\mu_0$  to the solution $w^1,\gamma^1,\mu^1$ of the system
\begin{equation}
\label{flohr:phiBRSF}
\left\{
\begin{aligned}
w^{1} &= w^0 + \dt \RS^\dx (w^0) + \dt\mu^1\partial_0w^0,\\
\mu^1 &= - \frac{\partial_0 \hat v^\top \big( w^0+\dt \RS^\dx (w^0) - \hat v \big)}{\dt \partial_0\hat v^\top\partial_0w},\\
\gamma^1 &= \gamma_0 + \dt\mu^1,
\end{aligned}\right. \qquad
\begin{aligned}
w^0 &= w_0,\\
 \gamma^0&=\gamma_0.
\end{aligned}
\end{equation}

Also for the sub-problem \eqref{flohr:A} we use a MOL approach, namely we spatially discretize \eqref{flohr:A} by finite differences,
 i.e. the discrete Laplacian with Dirichlet boundary conditions, $\partial_0^2$, is used to approximate the second spatial derivative,
 \begin{equation*}
 \frac{d}{dt}z_j = \partial_0^2 z_j, \quad z_j(0) = z^0_j.
 \end{equation*}
 For the time discretization we use backward Euler, because implicit methods have better stability properties for this type of equation.
 Using the linearity of $\partial_0^2$, this leads to $\phi_{A,\BE\dx}^\dt: (z_0,\gamma_0,\mu_0) \mapsto (z^1,\gamma^1,\mu^1)$ where
 \begin{equation}
 \left\{
 \begin{aligned}
  z^{1} &= (I-\dt\partial_0^2)^{-1}z^0,\\
 \gamma^1 &= \gamma^0,\\
 \mu^1 &= \mu^0,
 \end{aligned}\right. \qquad
 \begin{aligned}
 z^0 &= z_0,\\
  \gamma^0&=\gamma_0,\\
  \mu^0&=\mu_0,
 \end{aligned}
 \end{equation}
  such that $\phi_{A,\BE\dx}^\dt \approx \Phi^\dt_A$.

By using the Lie splitting \eqref{flohr:LieSplitting}, the full scheme for the freezing PDAE \eqref{flohr:freezedBurgers} is given by
\begin{equation}
\label{flohr:lo}
\tag{LO}
\begin{pmatrix}
v^{n+1}\\\gamma^{n+1}\\\mu^{n+1}
\end{pmatrix} \coloneqq \phi_{B,\RS\dx}^\dt \circ \phi_{A,\BE\dx}^\dt\begin{pmatrix}
v^{n}\\\gamma^{n}\\\mu^{n}
\end{pmatrix}
\end{equation}
for the orthogonal phase condition and by 
\begin{equation}
\label{flohr:lf}
\tag{LF}
\begin{pmatrix}
v^{n+1}\\\gamma^{n+1}\\\mu^{n+1}
\end{pmatrix} \coloneqq \varphi_{B,\RS\dx}^\dt \circ \phi_{A,\BE\dx}^\dt\begin{pmatrix}
v^{n}\\\gamma^{n}\\\mu^{n}
\end{pmatrix}
\end{equation}
for the fixed phase condition.

\subsection{Second order scheme.}
To construct a scheme with quadratic convergence in time and space we have to replace our numerical solution operators by suitable second order schemes and use Strang splitting instead of Lie splitting.
For the nonlinear hyperbolic part we use the second order semi-discrete scheme from \cite{kurganovTadmor}.
It is given as
\begin{equation}
\begin{aligned}
\frac{d}{dt}u_j(t) &=- \frac{1}{2\dx}\Big( f\big(u^+_{j+\slfrac{1}{2}}(t)\big) + f\big(u^-_{j+\slfrac{1}{2}}(t)\big) 
-  f\big(u^+_{j-\slfrac{1}{2}}(t)\big) - f\big(u^-_{j-\slfrac{1}{2}}(t)\big)\Big)\\
&\qquad+\frac{\kappa}{2\dx}\Big( u^+_{j+\slfrac{1}{2}}(t) - u^-_{j+\slfrac{1}{2}}(t)  - u^+_{j-\slfrac{1}{2}}(t) + u^-_{j-\slfrac{1}{2}}(t)\Big)\\
&\eqqcolon \KT^\dx(u(t)),
\end{aligned}
\end{equation}
where
\begin{equation*}
\begin{aligned}
u^\pm_{j+\frac{1}{2}}(t) &\coloneqq u_{j+\frac12 \pm \frac12}(t) \mp \frac{\dx}{2} (u_x)_{j+\frac12 \pm \frac12}(t)
\end{aligned}
\end{equation*}
for $j=-M,\dots,M$ with $u(t) \in \R^{2M+1}$ and $u_j(t) \in \R$ its $j$-th element.
The slopes are approximated using the minmod limiter
\[
(u_x)_j^n = \mm\left( \frac{u^n_j-u^n_{j-1}}{\dx}, \frac{u^n_{j+1}-u^n_j}{\dx} \right),
\]
where $\mm(a,b) \coloneqq \tfrac{1}{2}  [\sgn(a) + \sgn(b)] \cdot \min(|a|, |b|)$.
For the time integration we use Heun's method.
In the case of \eqref{pco},
$\phi_{B,\KT\dx}^\dt$ is the mapping $\phi_{B,\KT\dx}^\dt: (w_0,\gamma_0,\mu_0) \mapsto (w^1,\gamma^1,\mu^1)$
given by the solution of
\begin{equation}
\label{flohr:phiBKTF}
\left\{
\begin{aligned}
w^{*} &= w^0 + \dt \KT^\dx (w^0) + \dt\mu^*\partial_0w^0,\\
w^1 &= \tfrac{1}{2} w^0 + \tfrac{1}{2} \big( w^* + \dt \KT^\dx(w^*) + \dt \mu^* \partial_0w^* \big),\\
\mu^* &= - \frac{\partial_0 w^{0\top} \big( \partial_0^2 w^0 - w^0\partial_0 w^0 \big)}{\partial_0w^{0\top}\partial_0w^0},\\
\gamma^1 &= \gamma_0 + \dt\mu^1,\\
\mu^1 &= - \frac{\partial_0 w^{1\top} \big( \partial_0^2 w^1 - w^1\partial_0 w^1 \big)}{\partial_0w^{1\top}\partial_0w^1},
\end{aligned}\right. \qquad
\begin{aligned}
w^0 &= w_0,\\
 \gamma^0&=\gamma_0.
\end{aligned}
\end{equation}
For the fixed phase condition \eqref{pcf} we define $\varphi_{B,\KT\dx}^\dt$ 
as the mapping $(w_0,\gamma_0,\mu_0) \mapsto (w^1,\gamma^1,\mu^1)$
\begin{equation}
\label{flohr:phiBKTO}
\left\{
\begin{aligned}
w^{*} &= w^0 + \dt \KT^\dx (w^0) + \dt\mu^1\partial_0w^0,\\
w^1 &= \tfrac{1}{2} w^0 + \tfrac{1}{2} \big( w^* + \dt \KT^\dx(w^*) + \dt \mu^1 \partial_0w^* \big),\\
\mu^1 &= - \frac{\partial_0 \hat v^\top \big( w^0 + \dt \KT^\dx(w^0) - \hat v \big)}{\partial_0\hat v^\top\partial_0w},\\
\gamma^1 &= \gamma_0 + \dt\mu^1,
\end{aligned}\right. \qquad
\begin{aligned}
w^0 &= w_0,\\
 \gamma^0&=\gamma_0.
\end{aligned}
\end{equation}
 For the heat equation, we use the Crank-Nicolson\footnote{The Crank-Nicolson method used here is only the discretization in time by combining the forward and backward Euler method.} method to discretize in time and, as in the first order version, the discrete Laplacian with Dirichlet boundary conditions, $\partial_0^2$, is used in space.
The solution operator $\phi_{A,\CN\dx}^\dt$ is given by the mapping $(z_0,\gamma_0,\mu_0) \mapsto (z^1,\gamma^1,\mu^1)$
of 
 \begin{equation*}
 \left\{
 \begin{aligned}
  z^{1} &= (I-\tfrac{\dt}{2}\partial_0^2)^{-1}(I+\tfrac{\dt}{2}\partial_0^2)z^0,\\
 \gamma^1 &= \gamma^0,\\
 \mu^1 &= \mu^0,
 \end{aligned}\right. \qquad
 \begin{aligned}
 z^0 &= z_0,\\
  \gamma^0&=\gamma_0,\\
  \mu^0 &= \mu_0.
 \end{aligned}
 \end{equation*}
These methods where chosen, because they offer quadratic convergence 
for the individual problems
 and thus we can hope for quadratic convergence of the full problem with Strang splitting. 
Strang splitting \eqref{flohr:StrangSplitting} leads to our second order scheme given by
\begin{equation}
\label{flohr:so}
\tag{SO}
\begin{pmatrix}
v^{n+1}\\\gamma^{n+1}\\\mu^{n+1}
\end{pmatrix} = \phi_{B,\KT\dx}^\slfrac{\dt}{2} \circ \phi_{A,\CN\dx}^\dt \circ \phi_{B,\KT\dx}^\slfrac{\dt}{2}\begin{pmatrix}
v^{n}\\\gamma^{n}\\\mu^{n}
\end{pmatrix}
\end{equation}
for the orthogonal phase condition and by 
\begin{equation}
\label{flohr:sf}
\tag{SF}
\begin{pmatrix}
v^{n+1}\\\gamma^{n+1}\\\mu^{n+1}
\end{pmatrix} = \varphi_{B,\KT\dx}^\slfrac{\dt}{2} \circ \phi_{A,\CN\dx}^\dt \circ \varphi_{B,\KT\dx}^\slfrac{\dt}{2}\begin{pmatrix}
v^{n}\\\gamma^{n}\\\mu^{n}
\end{pmatrix}
\end{equation}
for the fixed phase condition.

\section{Numerical Results}

\begin{figure}[tb]
\centering

\renewcommand{\figureheight}{5.0cm}
\renewcommand{\figurewidth}{9.0cm}
\renewcommand{\yticklabelstyle}{yticklabel style={text width=0.65cm,align=right}}
\renewcommand{\xticklabelstyle}{} 

\input{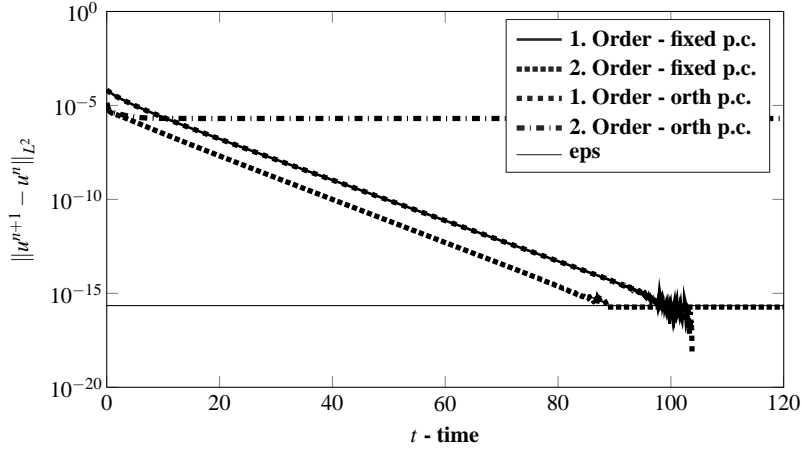}
\caption{Convergence to a numerical steady state except for the second order scheme with orthogonal phase condition.
}
\label{flohr:convtofp}
\end{figure}
The purpose of our schemes is to calculate viscous profiles by a simple forward simulation and thus we are interested in the quality of those profiles
obtained at the end of a long-time simulation. Note that we do not consider the initial convergence order on finite intervals.
For all following computations we use
the finite interval $[-15,15]$ and
 Dirichlet boundary conditions.
The reference function is given by \eqref{flohr:twSolutionBurgers} using $b = 1.5$, $c = -0.5$,
which is also used as initial value with $t=0$. This leads to a speed of $0.5$ for the traveling wave.
Since we are looking for numerical steady states in the co-moving frame,
we have to check if our numerical schemes yield steady states.
A steady state has the property $\frac{d}{dt}u(t) = 0$, which translates in the numerical case to $u^{n+1} = u^n$.
In Fig.~\ref{flohr:convtofp} we plot the time against the discrete $L^2$ distance $\|u^{n+1} - u^n\|_{L^2}$
and see that our schemes yield steady states at around $t \approx 100$ for \eqref{flohr:lo}, \eqref{flohr:lf} and \eqref{flohr:sf}
since $\|u^{n+1} - u^n\|_{L^2} \approx \text{machine precision}$.
For the Strang splitting scheme with orthogonal phase condition \eqref{flohr:so} we see that 
$\|u^{n+1} - u^n\|_{L^2}$ does not converge to zero and
the scheme does not offer a steady state. Solutions for this scheme leave the co-moving frame because the approximation of the speed is incorrect in this case.
For these computations we 
use $300$ steps in space and $\dt = \frac{\dx}{10}$.
 
\begin{figure}[tb]
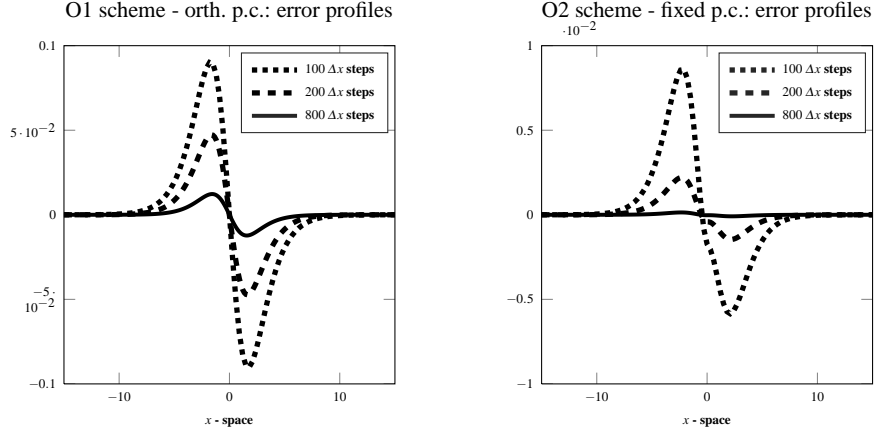

\renewcommand{\figureheight}{4.5cm}
\renewcommand{\figurewidth}{4.4cm}
\renewcommand{\yticklabelstyle}{yticklabel style={text width=0.65cm,align=right}}
\renewcommand{\xticklabelstyle}{} 

\input{figures/Figure_Poster_errorChangeO1.tikz} \hfill \input{figures/Figure_Poster_errorChangeO2.tikz}\\
\caption{Different numerical steady states for different $dt = \frac{\dx}{10}$. 
Note that the errors dominate where the profile varies the most and not at the boundary.}
\label{flohr:errorProfiles}
\end{figure}

Next, we consider the error profiles of the calculated steady states with different step sizes.
This result is shown in Fig.~\ref{flohr:errorProfiles}.
Obviously, we get different numerical steady states for different $dt = \frac{\dx}{10}$, which approximates the exact steady state better for smaller steps sizes.
In addition, we observe that the dominant error occurs in the profile and there is more-or-less no error at the boundary.

\begin{figure}[tb]
\centering

\renewcommand{\figureheight}{5.0cm}
\renewcommand{\figurewidth}{9.0cm}
\renewcommand{\yticklabelstyle}{yticklabel style={text width=0.65cm,align=right}}
\renewcommand{\xticklabelstyle}{} 

%
%
\definecolor{mycolor1}{rgb}{0.00000,0.44700,0.74100}%
\definecolor{mycolor2}{rgb}{0.85000,0.32500,0.09800}%
\definecolor{mycolor3}{rgb}{0.92900,0.69400,0.12500}%
\definecolor{mycolor4}{rgb}{0.49400,0.18400,0.55600}%
\begin{tikzpicture}

\begin{axis}[%
width=0\figurewidth,
height=0\figureheight,
at={(0\figurewidth,0\figureheight)},
scale only axis,
xmode=log,
xmin=0.0057692307692303,
xmax=0.300000000000001,
xminorticks=true,
xlabel={$10 \Delta t = \Delta x$ - step sizes},
ymode=log,
ymin=1e-05,
ymax=1,
yminorticks=true,
axis background/.style={fill=white},
title style={font=\bfseries},
title={discrete $L^2$-error to exact solution},
legend style={at={(0.03,0.97)},anchor=north west,legend cell align=left,align=left,draw=white!15!black},
title style={font=\small},xlabel style={font={\small\bfseries}},ylabel style={font=\small\bfseries},legend style={font=\tiny\bfseries},\yticklabelstyle,\xticklabelstyle,ticklabel style={font=\small}
]
\addplot [loosely dotted, every mark/.append style={solid, fill=gray, scale=1.3}, mark=triangle*]
  table[row sep=crcr]{%
0.211267605633802	0.140345268867059\\
0.149253731343283	0.100098416041098\\
0.10600706713781	0.0715725667991438\\
0.0748129675810478	0.0507574585657065\\
0.053003533568905	0.0360836075929751\\
0.0374531835205989	0.0255595774871515\\
0.0265017667844525	0.0181170561113963\\
0.0187382885696437	0.0128254802340609\\
0.0132567388422444	0.00908145940316351\\
0.00937207122774097	0.00642422330735846\\
0.00662836942112222	0.00454548050117969\\
0.00468676769254728	0.00321499027329723\\
0.00331455087835586	0.002274180677646\\
0.00234356690883608	0.00160821535941886\\
};
\addlegendentry{O1 scheme - orth. p.c.};

\addplot [dotted, every mark/.append style={solid, fill=gray, scale=1.2}, mark=square*]
  table[row sep=crcr]{%
0.211267605633802	0.00799778519934653\\
0.149253731343283	0.00400528394978266\\
0.10600706713781	0.00202114528487037\\
0.0748129675810478	0.00100661984762908\\
0.053003533568905	0.000505123358154324\\
0.0374531835205989	0.000252210598125602\\
0.0265017667844525	0.000126259599925078\\
0.0187382885696437	6.31171318977421e-05\\
0.0132567388422444	3.15884588870637e-05\\
0.00937207122774097	1.57870882925132e-05\\
0.00662836942112222	7.89642368378953e-06\\
0.00468676769254728	3.94783866376554e-06\\
0.00331455087835586	1.97474447854526e-06\\
0.00234356690883608	9.87262617074434e-07\\
};
\addlegendentry{O2 scheme - fixed p.c.};

\addplot [color=black,dashed,line width=1.0pt]
  table[row sep=crcr]{%
0.300000000000001	0.144\\
0.15	0.0720000000000002\\
0.0999999999999996	0.0479999999999998\\
0.0749999999999993	0.0359999999999997\\
0.0600000000000005	0.0288000000000002\\
0.0500000000000007	0.0240000000000003\\
0.0428571428571427	0.0205714285714285\\
0.0374999999999996	0.0179999999999998\\
0.0333333333333332	0.0159999999999999\\
0.0299999999999994	0.0143999999999997\\
0.0272727272727273	0.0130909090909091\\
0.0250000000000004	0.0120000000000002\\
0.023076923076923	0.011076923076923\\
0.0214285714285722	0.0102857142857147\\
0.0199999999999996	0.00959999999999979\\
0.0187500000000007	0.00900000000000034\\
0.0176470588235293	0.00847058823529409\\
0.0166666666666675	0.0080000000000004\\
0.0157894736842099	0.00757894736842076\\
0.0150000000000006	0.00720000000000027\\
0.0142857142857142	0.00685714285714283\\
0.0136363636363637	0.00654545454545456\\
0.0130434782608688	0.00626086956521704\\
0.0124999999999993	0.00599999999999966\\
0.0120000000000005	0.00576000000000022\\
0.0115384615384624	0.00553846153846194\\
0.0111111111111111	0.00533333333333331\\
0.0107142857142861	0.00514285714285734\\
0.0103448275862075	0.00496551724137959\\
0.00999999999999979	0.0047999999999999\\
0.00967741935483879	0.00464516129032262\\
0.00937500000000036	0.00450000000000017\\
0.0090909090909097	0.00436363636363666\\
0.00882352941176556	0.00423529411764747\\
0.0085714285714289	0.00411428571428587\\
0.00833333333333286	0.00399999999999977\\
0.00810810810810736	0.00389189189189153\\
0.00789473684210584	0.0037894736842108\\
0.00769230769230766	0.00369230769230768\\
0.00750000000000028	0.00360000000000014\\
0.00731707317073216	0.00351219512195144\\
0.00714285714285801	0.00342857142857184\\
0.00697674418604599	0.00334883720930208\\
0.00681818181818095	0.00327272727272685\\
0.00666666666666593	0.00319999999999965\\
0.00652173913043441	0.00313043478260852\\
0.00638297872340488	0.00306382978723434\\
0.00624999999999964	0.00299999999999983\\
0.00612244897959258	0.00293877551020444\\
0.00600000000000023	0.00288000000000011\\
0.00588235294117645	0.0028235294117647\\
0.0057692307692303	0.00276923076923055\\
};
\addlegendentry{linear ref.};

\addplot [color=black,dashdotted,line width=1.0pt]
  table[row sep=crcr]{%
0.300000000000001	0.00900000000000004\\
0.15	0.00225000000000001\\
0.0999999999999996	0.000999999999999993\\
0.0749999999999993	0.000562499999999989\\
0.0600000000000005	0.000360000000000006\\
0.0500000000000007	0.000250000000000007\\
0.0428571428571427	0.000183673469387754\\
0.0374999999999996	0.000140624999999997\\
0.0333333333333332	0.00011111111111111\\
0.0299999999999994	8.99999999999962e-05\\
0.0272727272727273	7.43801652892566e-05\\
0.0250000000000004	6.25000000000018e-05\\
0.023076923076923	5.32544378698221e-05\\
0.0214285714285722	4.59183673469423e-05\\
0.0199999999999996	3.99999999999983e-05\\
0.0187500000000007	3.51562500000027e-05\\
0.0176470588235293	3.11418685121105e-05\\
0.0166666666666675	2.77777777777805e-05\\
0.0157894736842099	2.49307479224357e-05\\
0.0150000000000006	2.25000000000017e-05\\
0.0142857142857142	2.0408163265306e-05\\
0.0136363636363637	1.85950413223141e-05\\
0.0130434782608688	1.70132325141758e-05\\
0.0124999999999993	1.56249999999982e-05\\
0.0120000000000005	1.44000000000011e-05\\
0.0115384615384624	1.33136094674576e-05\\
0.0111111111111111	1.23456790123456e-05\\
0.0107142857142861	1.14795918367356e-05\\
0.0103448275862075	1.07015457788359e-05\\
0.00999999999999979	9.99999999999957e-06\\
0.00967741935483879	9.36524453694084e-06\\
0.00937500000000036	8.78906250000067e-06\\
0.0090909090909097	8.26446280991847e-06\\
0.00882352941176556	7.78546712802919e-06\\
0.0085714285714289	7.34693877551076e-06\\
0.00833333333333286	6.94444444444365e-06\\
0.00810810810810736	6.57414170927563e-06\\
0.00789473684210584	6.23268698061033e-06\\
0.00769230769230766	5.91715976331357e-06\\
0.00750000000000028	5.62500000000043e-06\\
0.00731707317073216	5.35395597858484e-06\\
0.00714285714285801	5.10204081632776e-06\\
0.00697674418604599	4.86749594375265e-06\\
0.00681818181818095	4.64876033057732e-06\\
0.00666666666666593	4.44444444444347e-06\\
0.00652173913043441	4.25330812854394e-06\\
0.00638297872340488	4.07424173834393e-06\\
0.00624999999999964	3.90624999999956e-06\\
0.00612244897959258	3.74843815077142e-06\\
0.00600000000000023	3.60000000000027e-06\\
0.00588235294117645	3.46020761245672e-06\\
0.0057692307692303	3.32840236686337e-06\\
};
\addlegendentry{quad. ref.};

\end{axis}
\end{tikzpicture}%
\caption{Convergence rates of the numerical steady states to the exact one.
The scheme \eqref{flohr:lf} was omitted because it produces the same results as \eqref{flohr:lo}, whereas
the scheme \eqref{flohr:so} was ignored because it does not offer steady states.}
\label{flohr:convRates}
\end{figure}
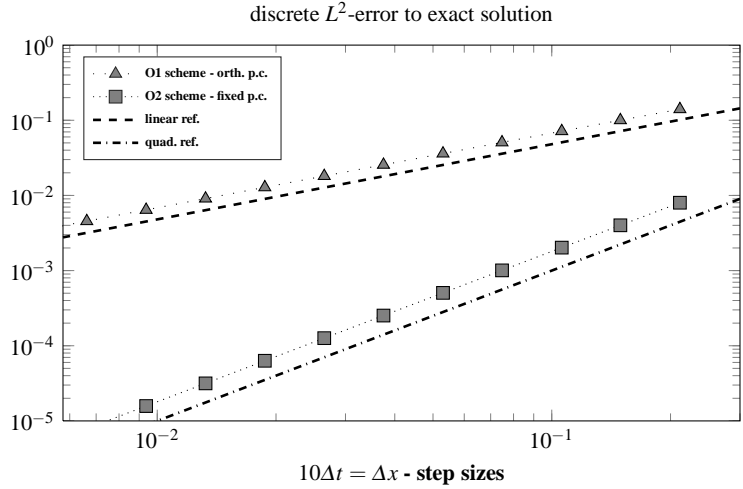
\begin{figure}[thb]
\centering

\renewcommand{\figureheight}{5.0cm}
\renewcommand{\figurewidth}{9.0cm}
\renewcommand{\yticklabelstyle}{yticklabel style={text width=0.65cm,align=right}}
\renewcommand{\xticklabelstyle}{} 

\input{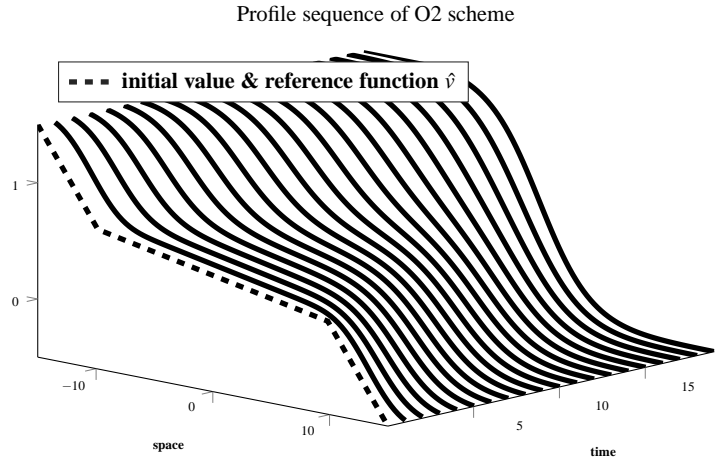}
\caption{Solution using initial value and reference function which only covers the rough behavior of the solution.}
\label{flohr:chooseRefFunc}
\end{figure}

The most interesting observation in our case is the convergence of our numerical steady states to the exact one.
For this we plot the discrete $L^2$ error of our states to the exact one for different step sizes in Fig.~\ref{flohr:convRates}.
Here we see linear convergence of our numerical steady states to the exact one for our first order scheme.
For the second order version we see quadratic convergence.

Finally, we note that usually
 the exact solution of the traveling wave is unknown.
 Therefore one has to guess some suitable reference function. In Fig.~\ref{flohr:chooseRefFunc} we see that a rough guess is sufficient for the initial value as well as for the reference function $\hat v$. The forward simulation approximates the traveling wave as before.

\begin{acknowledgement}
We gratefully acknowledge financial support by the Deutsche Forschungsgemeinschaft (DFG) through CRC 1173.
\end{acknowledgement}
%

%
%

\begin{thebibliography}{99.}%

    \bibitem{beynOttenRM2014}W.-J. Beyn and D. Otten and J.
    Rottmann-Matthes:
    \newblock{Stability and computation of dynamic patterns in {PDE}s.}
    \newblock{Current challenges in stability issues for numerical
    differential equations, 89--172}, {2014}.


    \bibitem{BT04}W.-J. Beyn and V. Th\"ummler:
    \newblock{Freezing solutions of equivariant evolution equations.}
    \newblock{SIAM J. Appl. Dyn. Syst., 3(2):85-116}, {2004}.

    \bibitem{F13}R. Flohr:
    \newblock{Splitting-Verfahren f\"ur partielle
    Differentialgleichungen mit Burgers-Nichtlinearit\"at.}
    \newblock{Masters Thesis, Bielefeld University}, {2013}.


    \bibitem{HoldenLubichRisebro}
    H. Holden and Chr. Lubich and N.H. Risebro:
    \newblock {Operator splitting for partial differential equations with Burgers nonlinearity.}
    \newblock {Math. Comput., 82:173-185}, {2013}.


    \bibitem{kurganovTadmor}A. Kurganov and E. Tadmor:
    \newblock {New High-Resolution Central Schemes for Nonlinear
    Conservation Laws and Convection-Diffusion Equations.}
    \newblock {J. Comput. Phys., 160:241-282}, {2000}.


  \bibitem{jens2012}J. Rottmann-Matthes:
  \newblock {Stability of parabolic-hyperbolic traveling waves.}
  \newblock {Dyn. Partial Differ. Equ., 9:29-62}, {2012}.
  
  
    \bibitem{jens2016}J. Rottmann-Matthes:
      \newblock {Freezing similarity solutions in multidimensional
      Burgers' equation.}
      \newblock {Preprint}, {2016}.
  \url{http://www.waves.kit.edu/downloads/CRC1173\_Preprint\_2016-27.pdf.},
  {2016}.
  
  
  \bibitem{stuart1996}A.M. Stuart and A.R. Humphries:
  \newblock {Dynamical systems and numerical analysis}
  \newblock {Cambridge Univ. Press}, {1996}.
  
  
  
  
  \bibitem{thuemmlerDiss}V. Th\"ummler:
          \newblock {Numerical analysis of the method of freezing
          traveling waves.}
          \newblock {Dissertation, Bielefeld University}, {2005}.

\end{thebibliography}
%
\biblstarthook{}

\end{document}